\documentclass{amsart}

\usepackage{amssymb}
\usepackage[all]{xy}

\newtheorem{thm}{Theorem}%[section]

\newtheorem{lem}[thm]{Lemma}
\newtheorem{cor}[thm]{Corollary}

\newtheorem{prop}[thm]{Proposition}

\newtheorem{conj}[thm]{Conjecture}
\theoremstyle{definition}
\newtheorem{defn}[thm]{Definition}

\newtheorem{say}[thm]{}
\newtheorem{exmp}[thm]{Example}

\newtheorem{prob}[thm]{Problem}
\newtheorem{ques}[thm]{Question}    %!!!!!!!!!!!!!!!!!!!!

\newtheorem{rem}[thm]{Remark}          

\newtheorem*{ack}{Acknowledgments}      % \renewcommand{\theack}{} 

\newtheorem{defn-thm}[thm]{Definition--Theorem}  %!!!!!!!!!!!!!!!!!!!!!!!!
\newtheorem{defn-lem}[thm]{Definition--Lemma}  %!!!!!!!!!!!!!!!!!!!!!!!!
\theoremstyle{remark}

\renewcommand{\c}[0]{{\mathbb C}}  

\renewcommand{\o}[0]{{\mathcal O}} 
\newcommand{\z}[0]{{\mathbb Z}}

\renewcommand{\r}[0]{{\mathbb R}} 

\renewcommand{\a}[0]{{\mathbb A}}

\newcommand{\p}[0]{{\mathbb P}}

\newcommand{\q}[0]{{\mathbb Q}}
\newcommand{\map}[0]{\dasharrow}
\newcommand{\qtq}[1]{\quad\mbox{#1}\quad}

\newcommand{\pic}[0]{\operatorname{Pic}}

\newcommand{\rank}[0]{\operatorname{rank}}

\newcommand{\supp}[0]{\operatorname{Supp}}    
\newcommand{\red}[0]{\operatorname{red}}    
    
\newcommand{\im}[0]{\operatorname{im}}

\newcommand{\sing}[0]{\operatorname{Sing}}

\newcommand{\chr}[0]{\operatorname{char}}

\newcommand{\hilb}[0]{\operatorname{Hilb}}

\newcommand{\len}[0]{\operatorname{length}}

\newcommand{\ord}[0]{\operatorname{ord}}

\newcommand{\tsum}[0]{\textstyle{\sum}}

\newcommand{\mtd}[1]{\mu(\operatorname{Td}_{#1})}

\def\loccoh#1.#2.#3.#4.{H^{#1}_{#2}(#3,#4)}

\DeclareMathAlphabet{\mathchanc}{OT1}{pzc}%
                                {m}{it}

\usepackage[all]{xy}\xyoption{dvips}

\newcommand{\elw}[0]{\operatorname{elw}}
\newcommand{\ndiv}{\hspace{-4pt}\not|\hspace{2pt}} 

\begin{document}
\bibliographystyle{amsalpha}

%\today

\title{Esnault-Levine-Wittenberg indices}
\author{J\'anos Koll\'ar}

\maketitle

These are notes of a lecture about the papers \cite{elw, wit}
and related results of \cite{mer-deg, ros-deg, MR2036596, MR2587339, MR3034413}.
The papers do not define {\it Esnault-Levine-Wittenberg indices,}
they focus on the two extreme cases
$\elw_0(X)$ (traditionally called the {\it index}) and  $\elw_{\dim X}(X)$.
In retrospect, the method of \cite[Sec.5]{MR2036596}
is equivalent to the computation of $\elw_1$ for Del~Pezzo surfaces.

The basic properties are at least implicitly in the above papers,
with the possible exception of the birational invariance
(\ref{elw.bir.lem}) and the  degree formula  (\ref{rost.prop}).
However, the definition and systematic use of the
ELW-indices streamlines several of the arguments.

It is interesting that the proof of  the birational invariance
does not rely on resolution. Instead, it uses what I call
the {\it Nishimura--Szab\'o lemmas}  (\ref{n-sz.lem}--\ref{k-sz.lem}).

 \begin{ack}
H.~Esnault, M.~Levine and O.~Wittenberg answered  many of my questions and
suggested several improvements. I tried to keep track and attribute
specific comments whenever possible.
Partial financial support   was provided  by  the NSF under grant number 
DMS-0968337.
\end{ack}

%It would be good  to understand if there is a closer relationship between
%the original Rost  degree formula  and the present variant.
 
\section{Definition and basic properties}

\begin{defn} \label{elw.defn}
Let $X$ be a proper scheme defined over a field $k$.
For $0\leq i\leq \dim X$ 
the {\it Esnault-Levine-Wittenberg index} $\elw_i(X) $
is defined as
$$
\elw_i(X):=\bigl(\chi(X, F): \dim F\leq i\bigr)\subset \z,
\eqno{(\ref{elw.defn}.1)}
$$
where  $F$ runs through all coherent sheaves  of dimension
$\dim F:=\dim\supp F\leq i$ over $X$. 
I will think of $\elw_i(X) $ as an ideal in $\z$; it can
be identified with its positive generator.
It is convenient to set $\elw_{(-1)}(X)=(0) $ and $ \elw(X):=\elw_{\dim X}(X)$.
It is clear that
$$
\elw_0(X)\subset \elw_1(X)\subset\cdots\subset \elw_{\dim X}(X).
\eqno{(\ref{elw.defn}.2)}
$$

If $X$ has a $k$-point $p\in X(k)$ then
$\chi\bigl(X, k(p)\bigr)=1$, hence
$\elw_0(X)= \elw_1(X)=\cdots= \elw_{\dim X}(X)=\z$.
Thus these notions are interesting only if $X$
does not have (or is not known to have) a $k$-point.

Usually $\operatorname{ind}(X):=\elw_0(X) $ is called the {\it index} of $X$.
Its generator is the smallest positive degree of a
(not necessarily effective) 0-cycle on $X$.

We see in (\ref{deviss.lem}.1) that if $X$ is integral then 
$$
\elw(X)=\bigl(\elw_{\dim X-1}(X),\chi(X, \o_X)\bigr).
\eqno{(\ref{elw.defn}.3)}
$$
This implies that if $\ell$ is a prime  such that
$\ord_{\ell}\elw(X)<\ord_{\ell}\elw_{\dim X-1}(X)$  then
 $$
\ord_{\ell}\elw(X)=\ord_{\ell}\chi(X, \o_X).
\eqno{(\ref{elw.defn}.4)}
$$
This is why the results of \cite{elw} involving $\chi(X, \o_X) $
are equivalent to statements about $\elw(X) $.
\end{defn}

\begin{rem} The above definition makes sense for a
proper scheme $X$ defined over a local Artin ring $A$. As a
consequence of (\ref{deviss.lem}.1) we see that
$\elw_i(X)=\elw_i(\red X)$ and $\red X$ 
can be viewed as a scheme over the residue field of $A$.
Once the basic results are established over Artin rings,
we concentrate on schemes over fields afterwards.
\end{rem}

The following   is  useful in computations.

\begin{lem} \label{chi.of.sheaf.cor} 
Let $X$ be a proper scheme over a local Artin ring $A$
and $F$  a coherent sheaf on $X$.
Let $Z_i\subset \supp F$  be the maximal dimensional irreducible components
with generic points $\eta_i$. Then
$$
\chi(X, F)-\tsum_i \len_{\eta_i}F\cdot\chi\bigl(Z_i, \o_{Z_i}\bigr)
\in \elw_{\dim  F-1}(X). 
\eqno{(\ref{chi.of.sheaf.cor}.1)}
$$
\end{lem}

Proof. The $K$-group of  coherent sheaves of  dimension
$\leq r$ is generated by the
$\o_Z$  where
$Z$ runs through all closed, integral subschemes of dimension
$\leq r$. On this group 
$$
F\mapsto \chi(X, F)-\tsum_i \len_{\eta_i}(F)\cdot\chi\bigl(Z_i, \o_{Z_i}\bigr)
\in \z/\elw_{r-1}(X)
$$
is linear. Thus it is enough to check that it vanishes on the generators.
If $\dim Z=r$ then 
$$
\tsum_i \len_{\eta_i}(\o_Z)\cdot\chi\bigl(Z_i, \o_{Z_i}\bigr)=\chi(X, \o_Z).
$$
If $\dim Z<r$ then $\tsum_i \len_{\eta_i}(F)\cdot\chi\bigl(Z_i, \o_{Z_i}\bigr)=0$
and $\chi(X, F)\in \elw_{r-1}(X)$. \qed

\begin{prop}[D\'evissage] \label{deviss.lem}
 Let $X$ be a proper scheme defined over a local Artin ring $A$.
The  ELW-indices
can also be computed the following ways.
\begin{enumerate}
\item $\elw_i(X) =\bigl(\chi(Z, \o_Z)\bigr)$ where
$Z$ runs through all integral subvarieties of dimension
$\leq i$.
\item $\elw_i(X) =\bigl(\chi(\bar Z, \o_{\bar Z})\bigr)$ where
$Z$ runs through all integral subvarieties of dimension
$\leq i$ and $\bar Z\to Z$ denotes the normalization.
\item $\elw_i(X) =\bigl(\chi( Z', \o_{Z'})\bigr)$ where
$Z$ runs through all integral subvarieties of dimension
$\leq i$ and $Z'\to Z$ is any proper birational morphism.
\item If $X$ is regular then  $\elw(X) =\bigl(\chi(X, E)\bigr)$ where
$E$ runs through all locally free sheaves on $X$.
\end{enumerate}
\end{prop}

Proof: The $K$-group of coherent sheaves of  dimension
$\leq i$ is generated by the
$\o_Z$  where
$Z$ runs through all integral subvarieties of dimension
$\leq i$. This implies (1).  If $X$ is regular then
locally free sheaves also generate the $K$-group, giving (4).

Finally (2--3) follow from this,
(\ref{elw.degd.lem}.3) and induction on $i$. \qed
\medskip

\begin{lem}\label{elw.degd.lem}
Let $f:X\to Y$ be a morphism of proper $k$-schemes
and $F$ a coherent sheaf on $X$.
Then
\begin{enumerate}
\item $\chi(X,F)\in \elw_{r}(Y)$ where $r=\dim f(\supp F)$.
\item If $Y$ is integral and $f$ is generically finite then \newline
$\chi(X,\o_X)-\deg (X/Y)\cdot \chi(Y,\o_Y)\in \elw_{\dim Y-1}(Y)$.
\item  If $f$ is birational then
$\chi(X,\o_X)- \chi(Y,\o_Y)\in \elw_{\dim Y-1}(Y)$.
\item  If $X,Y$ are normal and $f$ is birational then \newline
$\chi(X,\o_X)- \chi(Y,\o_Y)\in \elw_{\dim Y-2}(Y)$. 
\end{enumerate}
\end{lem}

Proof. By the Leray spectral sequence
$\chi(X,F)=\sum (-1)^i \chi(Y, R^if_*F)$
and the latter is in $\elw_{r}(Y)$ where
$r=\dim f(\supp F)$, showing (1).

If $X,Y$ are normal and $f$ is birational then
$\dim\supp R^if_*\o_X\leq \dim Y-2$ for $i>0$, thus
$\sum_{i\geq 1} (-1)^i \chi(Y, R^if_*\o_X)\in \elw_{\dim Y-2}(Y)$,
giving (4). 

 If $Y$ is integral and $f$ is generically finite then
the generic rank of $f_*\o_X$  equals the degree of $f$, thus
(2) follows from (\ref{chi.of.sheaf.cor}.1)
and (3) is a special case.\qed

\medskip

Applying (\ref{elw.degd.lem}.1) to all subvarieties we get the following.

\begin{cor} \label{elw.morph.lem}
Let $f:X\to Y$ be a morphism of proper $k$-schemes.
Then  $ \elw_i(X)\subset \elw_i(Y)$ for every $i$. \qed
\end{cor}

\begin{cor} Let $f:X\to Y$ be a  morphism of proper $k$-schemes
and $W\subset Y$ a closed subscheme such that $f$
 is an isomorphism over $Y\setminus W$. Then
$$
\elw_i(Y)=\elw_i(X)+\elw_i(W).
$$ 
\end{cor}

Proof. It is clear that $\elw_i(Y)\supset \elw_i(X)+\elw_i(W)$. 
Conversely, let $Z\subset Y$ be an integral subvariety.
If $Z\subset W$ then set $Z'=Z$ and note that
$\chi(Z', \o_{Z'})\in \elw_i(W)$. 

If $Z\not\subset W$ then let $Z'$ be the birational transform
of $Z$ in $X$. Thus $\chi(Z', \o_{Z'})\in \elw_i(X)$.
Using (\ref{elw.degd.lem}.3) we see that
$\elw_i(Y)\subset \elw_i(X)+\elw_i(W)$. \qed
\medskip

A more interesting variant of 
(\ref{elw.morph.lem}) is the following.

\begin{prop} \label{elw.map.lem} 
Let $f:X\map Y$ be a rational map of proper $k$-schemes.
Assume that $X$ is regular.
Then  $ \elw_i(X)\subset \elw_i(Y)$ for every $i$.
\end{prop}

Proof. By induction on $i$. Let $Z\subset X$ be an  integral subvariety
of dimension $i$. By (\ref{n-sz.lem}) 
we have a birational morphism $Z'\to Z$ and a
morphism $Z'\to Y$. Thus  $\chi(Z', \o_{Z'})\in \elw_i(Y)$ 
and 
$\chi(Z', \o_{Z'})-\chi(Z, \o_{Z})\in \elw_{i-1}(X)\subset \elw_{i-1}(Y)$
by induction. 
Thus $\chi(Z, \o_{Z})\in\elw_i(Y)$.\qed

\begin{cor} \label{elw.bir.lem} 
For every $i$, the ELW-index  $X\mapsto  \elw_i(X)$
is a birational invariant of proper, regular $k$-schemes. \qed
\end{cor}

An immediate consequence of 
(\ref{elw.degd.lem}.1) and  (\ref{elw.map.lem})
is the following quite useful result.

\begin{cor} Let $X$ be a proper, regular $k$-scheme and
$g:X\map Y$ a map. 
If $\chi(X, \o_X)\not\in \elw_r(Y)$ then
$\dim g(X)>r$. \qed
\end{cor}

For $i=\dim X$,
the following  degree formula is in
\cite{mer-deg, MR2587339, MR3034413}.
(See (\ref{todd.gives.0.cyc.lem.n-1}) for its precise relationship to the
 versions given in \cite{MR2587339, MR3034413}.)

\begin{prop}\label{rost.prop}
Let $f:X\map Y$ be a generically finite, 
rational map of proper $k$-schemes
of the same dimension.
Assume that $Y$ is integral and regular. Then
$$
\deg (X/Y)\cdot \elw_i(Y)\subset  \elw_i(X)+\elw_{i-1}(Y).
$$
\end{prop}

Proof.  Let $Z\subset Y$ be an  integral subvariety
of dimension $i$. By (\ref{k-sz.lem}) 
we have a generically finite  morphism $Z'\to Z$ of degree $d$ and a
morphism $Z'\to X$. Thus  $\chi(Z', \o_{Z'})\in \elw_i(X)$ 
by (\ref{elw.degd.lem}.1)
and 
$\chi(Z', \o_{Z'})-\deg (X/Y)\cdot \chi(Z, \o_{Z})\in \elw_{i-1}(Y)
\subset \elw_{i-1}(Y)$
by (\ref{elw.degd.lem}.2). 
Thus $\deg (X/Y)\cdot \chi(Z, \o_{Z})\in\elw_i(X)+\elw_{i-1}(Y)$.\qed

\begin{cor}\label{rost.cor}
Notation and assumptions as in (\ref{rost.prop}). 
Fix a prime $\ell$ such that
$\ord_{\ell}\elw_{i-1}(Y)>\ord_{\ell}\elw_{i}(Y)$ and
$\ell\ndiv \deg (X/Y)$.

Then $\ord_{\ell}\elw_{i}(X)=\ord_{\ell}\elw_{i}(Y)$. \qed
\end{cor}

\subsection*{Cycle class map}{\ }

The ELW-indices determine how far the Euler characteristic is from
being linear on the Chow groups.

\begin{defn}
For an algebraic cycle $Z=\sum_ia_iZ_i\subset X$ set
$$
\chi(Z):=\tsum_ia_i\chi(Z_i, \o_{Z_i})\qtq{and}
\chi(\bar Z):=\tsum_ia_i\chi(\bar Z_i, \o_{\bar Z_i}),
$$
where $\bar Z_i\to Z_i$ denotes the normalization.

If $W\subset X$ is an integral subvariety of dimension $r$ then 
 $\chi(\bar W, \o_{\bar W})-\chi(W, \o_{W})\in \elw_{r-1}(X)$
by (\ref{elw.degd.lem}.3). Thus
$$
\chi(\bar Z)-\chi(Z)\in \elw_{\dim Z-1}(X).
$$
\end{defn}

\begin{prop}\label{chi.on.alg.eq}
 Let  $B_r(X)$ denote the group of $r$ dimensional cycles in $X$
modulo algebraic equivalence. Then
$Z\mapsto \chi(Z)$ and  $Z\mapsto \chi(\bar Z)$ define the same
well-defined linear map
$$
\chi: B_r(X)\to  \elw_r(X)/\elw_{r-1}(X).
$$
\end{prop}

Proof. 
Two  $r$-cycles $Z^1, Z^2$ are algebraically equivalent
if there are  
\begin{enumerate}
\item an irreducible, nonsingular  curve $C$ with two
 points $p_1, p_2\in C(k)$,
\item a  flat and  proper morphism $g:W\to C$,
\item  a morphism  $\pi: W\to X$  and
\item an effective  $r$-cycle $Z^c$ such that 
$\pi_*\bigl[g^{-1}(p_i)\bigr]=Z^i+Z^c$ for $i=1,2$.
\end{enumerate}
Set $W^i:=g^{-1}(p_i)$ and
let $W^i_j\subset W^i$ be the irreducible components
with multiplicities  $m^i_j$. By (\ref{chi.of.sheaf.cor}.1)
$$
\chi\bigl( W^i, \o_{W^i}\bigr)-\tsum_j m^i_j\cdot\chi\bigl( W^i_j, \o_{W^i_j}\bigr)
\in \elw_{r-1}(W)
$$
for $i=1,2$. Furthermore, by (\ref{elw.degd.lem}.2)
$\chi\bigl( \pi_*[W^i_j]\bigr)-\chi\bigl( W^i_j\bigr)
\in \elw_{r-1}(X)$.
Therefore
$$
\chi\bigl( Z^i\bigr)+\chi\bigl( Z^c\bigr)-
\tsum_j m^i_j\cdot\chi\bigl( W^i_j, \o_{W^i_j}\bigr)
\in \elw_{r-1}(X).
$$
Finally  $\chi\bigl( W^1, \o_{W^1}\bigr)=\chi\bigl( W^2, \o_{W^2}\bigr)$
since $g$ is flat. \qed

\begin{exmp} Let $Q^n\subset \p^{n+1}$ be the (empty) quadric
$(x_0^2+\cdots+x_{n+1}^2=0)$ over $\r$. Note that  the quadric $Q^2$
 is isomorphic to the  product $Q^1\times Q^1$.
Let $C_2, C_4\subset Q^1\times Q^1$ be  rational curves of bidegrees (1,1)
(resp.\  (1,3)).

Both of these can be viewed as curves in $Q^3$ since
$Q^2\subset Q^3$.
The cycles $C_4$ and $2\cdot C_2$ have the same degree, hence
they are algebraically equivalent  over $\c$. 

However, they are  not algebraically equivalent over $\r$
since  $\chi(C_4)=1$ but $2\chi(C_2)=2$.
\end{exmp}

\begin{rem}[ML] This shows that 
$F\mapsto \chi(X,F)\in \z/\elw_{\dim F-1}(X)$ 
can be viewed as the image of the push-forward map from the
connective algebraic K-theory of $X$ to the base-field. This is given by
taking $CK^q(X)$ to be the image of $K_0(M^q(X))$ in  $K_0(M^{q-1}(X))$, where
$M^q(X)$ is the category of coherent sheaves on X with support in codim at
least q. This theory was first defined by \cite{MR2400737}. 
\cite{Dai-levine}   shows that $CK^*$ is the
universal theory with formal group law $u+v-buv \in \z[b][[u,v]]$. 
\end{rem}

\section{The Nishimura--Szab\'o lemmas on rational correspondences}{\ }

\cite{nis} proved that if a regular $k$-scheme has a $k$-point
then any proper $k$-scheme birational to it also  has a $k$-point.
Around 1992  Endre Szab\'o found a new short argument.
 (The proof is reproduced in \cite[Prop.A.6]{MR1782331-app}
and \cite[p.183]{ksc}.) 
Another application of the method was also used in \cite{MR1782331-app};
its generalization
(\ref{k-sz.lem}) is needed to prove the 
  degree formula  (\ref{rost.prop}).

\begin{lem} \label{n-sz.lem}
Let $f:X\map Y$ be a rational map of proper $k$-schemes.
Let $Z\subset X$ be a closed,  integral subscheme that is not contained
in $\sing X$. Then there is a birational morphism $Z'\to Z$ 
such that there exists a
morphism $f'_Z:Z'\to Y$.
\end{lem}

We make no claim about $f'_Z$ beyond its existence.
In particular, it can be a constant map. Thus (\ref{n-sz.lem})
is interesting only if $Y$ is not known to have a $k$-point.
\medskip

Proof. By induction on $\dim X$. If $\dim X\leq \dim Z+1$ then
$f$ is defined at the generic point of $Z$, hence we can take
$Z'$ to be the closure of the graph of $f|_Z$.

If $\dim X> \dim Z+1$, take the blow-up $B_ZX\to X$ and let
$E_Z\subset B_ZX$ be the unique irreducible component
of the exceptional divisor that dominates $Z$. Note that
$E_Z\to Z$ is generically a projective space bundle, hence
it has a rational section $Z_1\subset E_Z$ that  maps
birationally to $Z$. Note that $f$ gives a rational map
$f_E:E_Z\map Y$. Induction gives  
a birational morphism $Z'_1\to Z_1\to Z$ and a
morphism $Z'_1\to Y$. \qed

\medskip
The following is a variant of \cite{MR1782331-app}.

\begin{lem}%[Koll\'ar--Szab\'o lemma] 
\label{k-sz.lem}
Let $f:X\map Y$ be a generically finite, rational map of proper $k$-schemes
of the same dimension.
Assume that $Y$ is integral. 
Let $Z\subset Y$ be a closed,  integral subscheme that is not contained
in $\sing Y$. Then there are
\begin{enumerate}
\item  a reduced, proper $k$-scheme $Z'$,
\item a generically finite morphism $Z'\to Z$ such that
$\deg (Z'/Z)=\deg (X/Y)$ and
\item a morphism  $Z'\to X$.
\end{enumerate}
\end{lem}

Proof. By induction on $\dim Y$. 

Replacing $X$ by the normalization of the closure of the graph of $f$, 
we may assume that $X$ is normal and 
$f:X\to Y$ is a morphism.
 
If  $\dim Y\leq \dim Z+1$ then 
$f$ is finite over  the generic point of $Z$.
Let $Z_i\subset X$ be the irreducible components of $f^{-1}(Z)$
that dominate $Z$ and $e_i$ the ramification index of $f$ along $Z_i$.
Then  $\deg (X/Y)=\sum_i e_i\deg(Z_i/Z)$, thus we can take
$Z'$ to be the disjoint union of $e_i$ copies of $Z_i$ for every $i$.

If $\dim Y> \dim Z+1$, take the blow-up $B_ZY\to Y$ and let
$E_Z\subset B_ZY$ be the unique irreducible component
of the exceptional divisor that dominates $Z$. Note that
$E_Z\to Z$ is generically a projective space bundle, hence
it has a rational section $Z_1\subset E_Z$ that  maps
birationally to $Z$. 

Replace $X$ by the  by the normalization of the closure of the graph of 
$X\map B_ZY$. Let  $F_i\subset X$ be the irreducible components of $f^{-1}(E_Z)$
that dominate $E_Z$ and $e_i$ the ramification index of $f$ along $F_i$. 
By induction, there are
$Z'_{1i}\to Z_1$ and morphisms $ Z'_{1i}\to F_i\to X$. 
Thus  we can take
$Z'$ to be the disjoint union of $e_i$ copies of $Z'_{1i}$ for every $i$.
 \qed

\begin{rem} The computations of \cite[Lem.5.3]{wit} show that 
the above results also hold if the ambient variety  has quotient singularities
at the generic point of $Z$ and $k(Z)$ is perfect. 
(For $\dim Z>0$ this restricts us to characteristic 0.)

Here quotient singularity is understood in the strong sense: it should be
Zariski locally a quotient of a regular scheme. 
The analogous assertion is not true for singularities that are
quotients only \'etale locally. For instance, take 
$X=(x^2+y^2+z^2=0)\subset \a^3_{\r}$ and $Z=(0,0,0)$. After blowing up the
origin, we get a surface with no real points. On the other hand,
over $\c$ the singularity is isomorphic to $\c^2/(u,v)\sim(-u,-v)$.
\end{rem}

\section{Examples}

Proposition \ref{chi.on.alg.eq} makes it relatively easy to compute
the ELW-indices when generators of the groups  $B_r(X)$
are known.

\begin{exmp}[ML] Let $p$ be a prime and $X$  a nontrivial Severi--Brauer variety
of dimension $p-1$. Then
$$
\elw_{0}(X)=\cdots=\elw_{p-2}(X)=p\qtq{and} \elw_{p-1}(X)=1.
$$
To see this note that a Severi--Brauer variety of dimension
$n-1$ has an (effective) 0-cycle of degree $n$ and it has a 
0-cycle of degree 1 iff it is trivial. Thus 
$\elw_{0}(X)=p$. At the other end, $\chi(X,\o_X)=1$ shows that
$\elw_{p-1}(X)=1 $. Thus the only question is when $\elw_i$ drops from
$p$ to $1$. Assume that 
$\elw_{i-1}(X)=p $.

Let $K/k$ be a splitting field of degree $p$ and
$L_K\subset X_K$ a linear subspace of dimension $i$.
Let $L_1,\dots, L_p$ be its conjugates over $k$.
Then $Z:=L_1+\cdots + L_p$ is defined over $k$ and it generates
the Chow group $A_i(X)$ if $i<\dim X$. % by \cite{??}.  
$Z$ could be singular, but
its normalization $\bar Z$ has Euler characteristic
$\chi(\bar Z, \o_{\bar Z})=p\chi(L_K, \o_{L_K})=p$. 

By (\ref{chi.on.alg.eq}), 
$Z\mapsto \chi(\bar Z, \o_{\bar Z})\in \z/p$ is a well defined
linear map on  $A_i(X)$ that vanishes on the generator of $A_i(X)$.
Thus $\elw_{i}(X)=p $.

(OW) notes that this is also a direct consequence of
(\ref{todd.gives.0.cyc.lem.n-1}).
\end{exmp}

Next we compute the ELW-indices for certain products of 
 general curves.

%This needs a modified version of $\elw_0$.

%\begin{defn} Let  $C$ be a 1-dimensional proper scheme
%over a field $k$. Set
%$$
%\elw_0^*(C):=\bigl(m: \pic^m(C_{\bar k})(k)\neq \emptyset\bigr).
%$$
%It is clear that  $\elw_0^*(C)\supset \elw_0(C)$ and frequently the
%two are equal. However, if $C$ is smooth and geometrically rational then
%$\elw_0^*(C)=(1)$ yet $\elw_0(C)=(2)$ if $C(k)=\emptyset$.
%\end{defn}

\begin{prop}\label{prod.curves.prop}
  Let $k$ be a field and
$C_1,\dots, C_n$ smooth, irreducible, projective curves over $k$. Assume that
\begin{enumerate}
\item  $\elw_0(C_i)\subset (m)$ for every $i$ for some $m\geq 1$ and
\item  $\pic\bigl(C_1\times\cdots\times C_n)
=\pi_1^*\pic(C_1)+\cdots +\pi_n^*\pic(C_n)$ where  $\pi_i$ denotes the $i$-th
coordinate projection.
%the $C_1,\dots, C_n$ satisfy the genericity 
%conditions of (\ref{pic.of.prod.lem}).
\end{enumerate}
Then  $\elw_i\bigl(C_1\times\cdots\times C_n)\subset (m)$ for $i<n$.
\end{prop}

Proof. Set $X_r:=C_1\times\cdots\times C_r$. The proof is by induction on $r$.

Let $F$ be a coherent sheaf on $X_n$.
Consider the coordinate projection
$\Pi_n: X_n\to X_{n-1}$. 
If $\dim F\leq n-2$ then  $\dim \Pi_n(\supp F)\leq n-2$, hence,
 by (\ref{elw.degd.lem}.1) and induction,
$\chi(X_n, F)\in \elw_{n-2}\bigl(X_{n-1}\bigr)\subset (m)$.

We are left with the case when $\dim  F=n-1$.
By (\ref{chi.of.sheaf.cor}.1) it is enough to show that
$\chi(D, \o_D)\in (m)$ for every effective divisor $D\subset X_n$.
Using the exact sequence
$$
0\to \o_{X_n}(-D)\to \o_{X_n}\to \o_D\to 0
$$
we are reduced to proving that
$$
\chi(X_n, L)\equiv \chi(X_n, \o_{X_n})\mod m
$$
for every line bundle $L$ on $X_n$. By assumption (2),
there are line bundles $L_i$ on $C_i$ such that 
$L\cong \otimes_i\pi_i^*L_i$.
Therefore
$$
\chi(X_n, L)=\textstyle{\prod}_i \chi(C_i, L_i)=
\textstyle{\prod}_i\bigl(\chi(C_i, \o_{C_i})+\deg L_i\bigr).
$$
By assumption (1), $\deg L_i\in (m)$ for every $i$,  thus
$$
\chi(X_n, L)\equiv \textstyle{\prod}_i\ \chi(C_i, \o_{C_i})=
\chi(X_n, \o_{X_n})\mod m.\qed
$$
\medskip

In applications the tricky part is to check the condition
(\ref{prod.curves.prop}.2). Let us start over algebraically
closed  fields.

\begin{lem} \label{pic.of.prod.lem}
Let
$Y_i$ be normal, irreducible, proper varieties over an algebraically
closed  field $k$.
The following are equivalent.
\begin{enumerate}
\item $\pic\bigl(Y_1\times\cdots\times Y_n)
=\pi_1^*\pic(Y_1)+\cdots +\pi_n^*\pic(Y_n)$ where  $\pi_i$ denotes the $i$th
coordinate projection.
\item For $i\neq j$, every morphism $Y_i\to \pic(Y_j)$ is constant.
\item For $i\neq j$, every morphism $\pic^0(Y_i)\to \pic^0(Y_j)$ is constant.
\qed
\end{enumerate}
\end{lem}

The above conditions hold if the $\pic^0(Y_i)$ are sufficiently general
and independent, but they may be hard to check in concrete
situations. Consider the case when $Y_i=C_i$ are smooth curves.
For very general curves $\pic^0(C_i)$ is a simple Abelian variety;
see \cite{MR0453754, MR529440, MR1748293, MR2040573} for explicit examples over $\q$.
Hence if, in addition,  the $C_i$ all have different genera,
then (\ref{pic.of.prod.lem}.3) holds.
The same for sufficiently general
and independent curves of the same genus $>0$.

Over arbitrary fields, an extra complication comes from the
Brauer group.

\begin{say}[Map to the Brauer group]\label{brauer.map.say}
  Let
$Y$ be normal, irreducible, proper variety over a  perfect field $k$
with algebraic closure $\bar k$. Let $L$ be a line bundle on $Y_{\bar k} $
that is isomorphic to its Galois conjugates. Equivalently,
$[L]\in \pic\bigl(Y_{\bar k}\bigr)(k)$.
If the linear system $|L|$ is nonempty, then,  as a subscheme of $\hilb(Y)$;
it  is
defined over $k$ and is isomorphic to a projective space
over $\bar k$. Thus $|L|$ defines an element of the Brauer group
$ \operatorname{Br}(k) $.
This gives an exact sequence
$$
\pic(Y)\to \pic\bigl(Y_{\bar k}\bigr)(k)
\stackrel{\operatorname{br}_Y}{\longrightarrow} \operatorname{Br}(k).
$$
(See \cite[Chap.8]{blr} for a more conceptual construction.)
It is clear that if
$[L_i]\in \pic\bigl((Y_i)_{\bar k}\bigr)(k)$
then
$$
\operatorname{br}_{\prod Y_i}\bigl(\otimes_i\pi_i^*L_i\bigr)=
\tsum_i \operatorname{br}_{Y_i}\bigl(L_i).
$$
\end{say}

\begin{lem} \label{pic.of.prod.lem.2}
Let
$Y_1,\dots, Y_n$ be normal, irreducible, 
proper varieties over a  perfect field $k$.
Assume that they satisfy the equivalent conditions
(\ref{pic.of.prod.lem}.1--3) over $\bar k$. 
The following are equivalent.
\begin{enumerate}
\item $\pic\bigl(Y_1\times\cdots\times Y_n)
=\pi_1^*\pic(Y_1)+\cdots +\pi_n^*\pic(Y_n)$.
\item The subgroups
$\im\bigl[\operatorname{br}_{Y_i}: \pic\bigl((Y_i)_{\bar k}\bigr)(k)
\to \operatorname{Br}(k)\bigr]$
are independent  in $\operatorname{Br}(k)$.
(Subgroups  $A_i$ of an Abelian group $A$
are {\it independent} if  $a_i\in A_i$ and $\sum a_i=0$ implies that
$a_i=0$ for every $i$.)
\qed
\end{enumerate}
\end{lem}

\begin{exmp}[Products of conics] Let $C_i\subset \p^2$ be plane conics
over  a  perfect field $k$. Set $Y_n=C_1\times\cdots\times C_n$.
The image of $\operatorname{br}_{C_i} $
equals the subgroup generated by $C_i$ in $\operatorname{Br}(k)$.
Thus we see that the following are equivalent.
\begin{enumerate} 
\item $\elw_{n-1}(Y_n)\subset (2)$ and
\item the classes $[C_i]\in \operatorname{Br}(k)$ are
nonzero and independent.
\end{enumerate}
If we are over $\q$ then any finite collection of conics
has a point in a quadratic extension. Thus we obtain the following.

Set $C_i:=(x^2+y^2=p_iz^2)$ where the $p_i$ are distinct primes
congruent to $3$ modulo $4$ and   $Y_n:=C_1\times\cdots\times C_n$.
Then
$$
\elw_{0}(Y_n)=\cdots =\elw_{n-1}(Y_n)=(2)
\qtq{and} \elw_{n}(Y_n)=(1).
$$
\end{exmp}

\begin{exmp}[Products of hyperelliptic curves]
We work over the field $k=\c(t)$. 
A hyperelliptic curve over $k$ is given by an equation
$$
C:=\bigl(z^2=f(x,y)\bigr)\subset \p^2(1,1,m) 
$$
where $f(x,y)\in \c(t)[x,y]$ is homogeneous of degree $2m$.
We will look at it as a cover of $\p^1\times \p^1$ where
$\a^2$ (with coordinates $t,x/y$) is an affine chart.
Thus we have a double cover $\pi:S\to \p^1\times \p^1$
with branch curve $B\subset \p^1\times \p^1$.
For very general $S$,
$$
\pi^*\pic \bigl (\p^1\times \p^1\bigr)\to \pic(S)
\qtq{is an isomorphism;}
$$
see \cite{MR699163, rav-sri}. This implies that, for very general $f$,
$\pic(C)=\z\bigl[\pi^*\o_{\p^1}(1)\bigr]$. 
In particular, if $m$ is odd then 
$\elw_0(C)=2$ and $\elw_0(C)=1$.

For stronger examples over $\q$, see \cite{MR2221085, MR2578467}.

Note that $\pic(C)$ can be viewed as the generic fiber of
the family of the Picard varieties of the fibers of
$p_1\circ \pi:S\to \p^1$ where  $p_1$ is  the first
coordinate projection of $\p^1\times \p^1$. 
The singular fibers correspond to the
branch  points of $p_1:B\to \p^1$. 

Let us now choose our curves  $C_i$ such that the
branch  points of the corresponding  $p_1:B_i\to \p^1$ are all different.
Then the different $\pic^0(C_i)$ satisfy 
the condition (\ref{pic.of.prod.lem}.3). The Brauer group of
$\c(t)$ is trivial. Thus if 
 $Y_n:=C_1\times\cdots\times C_n$ is the product of such generic  curves
then
$$
\elw_{0}(Y_n)=\cdots =\elw_{n-1}(Y_n)=(2)
\qtq{and} \elw_{n}(Y_n)=(1).
$$
\end{exmp}

\begin{exmp}[Products of real curves] Let $C$ be
a geometrically irreducible smooth real curve. The interesting case
is when $C$ has  no real points and even genus. Equivalently, when
$\elw_0(C)=(2)$ and $\elw_1(C)=(1)$. Note that
$\operatorname{Br}(\r)=\z/2$ and the Brauer map
${\operatorname{br}_C}: \pic\bigl(C_{\c}\bigr)(\r)
\to  \operatorname{Br}(\r)=\z/2$
is surjective, essentially by a result that goes back to   \cite{witt-br}
(see also \cite{MR0412193}). 
Thus we see that for any product $Y_n$ of such curves we have
$$
\elw_{0}(Y_n)=(2) \qtq{and}\elw_{1}(Y_n)= \cdots =
 \elw_{n}(Y_n)=(1).
$$
\end{exmp}

\section{The sequence of ELW-indices}

There are several general questions 
about the ELW-indices that should be explored.
Esnault asks about all possible
 sequences
$\elw_0(X),\dots, \elw_{\dim X}(X)$.

\begin{defn}\label{mu.td.n}
Let $\mtd{n}$ denote  the denominator appearing in the
Todd class in dimension $n$. By \cite[1.7.3]{hirz}, it is
$$
\mtd{n}=\textstyle{\prod}_p\ p^{[n/(p-1)]}\qtq{where $p\leq n+1$ is a prime.}
\eqno{(\ref{mu.td.n}.1)}
$$
The sequence starts as
$\mtd{1}=2$, $\mtd{2}=12$, $\mtd{3}=24$, $\mtd{4}=720$.

Note that $\mtd{n}$ is very close to $n!$. Indeed,
we can rewrite the formula as
$$
\mtd{n}=\textstyle{\prod}_p\ p^{[\frac{n}{p}+\frac{n}{p^2}+\cdots]}
\qtq{while}
n!=\textstyle{\prod}_p\ p^{[\frac{n}{p}]+[\frac{n}{p^2}]+\cdots}.
$$
Each $\mtd{n}$ divides any later $\mtd{m}$.
There is also the more delicate relation
$$
n!\cdot \mtd{m}\ \big| \ \mtd{n+m-1}.
\eqno{(\ref{mu.td.n}.2)}
$$
\end{defn}

Without much evidence, let me propose the following.

\begin{prob} \label{HE.conj}
Let $e_0,\dots, e_n$ be a sequence of natural numbers.
Then there is a field $k$ and a $k$-scheme  (or smooth $k$-variety)
$X$ of dimension $n$ such that $\elw_r(X)=(e_r)$ for every $r$ iff
the following hold.
\begin{enumerate}
\item $e_{r+1}|e_r$ for every $r$ and
\item $e_0| \mtd{r}\cdot e_r$ for every $r$.
\end{enumerate}
\end{prob}

Next we discuss some evidence supporting the above formulation of 
(\ref{HE.conj}).
The necessity of $e_{r+1}|e_r$  follows from (\ref{elw.defn}.2).
The divisibility $e_0| \mtd{r}\cdot e_r$ is more subtle and
I do not know a complete proof.

\begin{lem} \label{todd.gives.0.cyc.lem} Let $X$ be a proper $k$-scheme over a
field of characteristic 0. Then 
 $$
\mtd{r}\cdot \elw_{r}(X) \subset \elw_0(X).
\eqno{(\ref{todd.gives.0.cyc.lem}.1)}
$$
\end{lem}

Proof.  If $W$ is smooth, proper and of dimension $r$
then Riemann--Roch says that  $\mtd{r}\cdot \chi(W, \o_W)$ is 
a polynomial in the Chern classes. In particular, $W$ has a
0-cycle of degree $\mtd{r}\cdot \chi(W, \o_W)$.

By (\ref{deviss.lem}.3), $\elw_{r}(X) $ is generated by the 
$\chi( Z', \o_{Z'})$ where
$Z$ runs through all integral subvarieties of dimension
$\leq r$ and $Z'\to Z$ is any resolution.\qed
\medskip

The only missing ingredient in positive characteristic is
resolution of singularities. The method of \cite{elw}
shows that 
one can use de~Jong's alterations \cite{deJ-alt, MR1450427}
(in the stronger form proved by Gabber, cf.\
  \cite[Exp.IX]{2012arXiv1207.3648I}) to 
prove  that (\ref{todd.gives.0.cyc.lem}.1)
holds in $\z[p^{-1}]$.
%prime-to-$p$ part in $\chr p$.  

\medskip

A more interesting part of (\ref{HE.conj}) is the claim
that there are no additional relations between the
$\elw_i$. I do not even have a plausible argument, but the following example
suggests that,  for $n=\dim X$, there should not be any relations between
$\elw_n(X)$ and $\elw_{n-1}(X)$ in general.

\begin{exmp}
Let $X$ be a smooth projective variety of dimension $n$ such that
$\pic(X)=\z[H]$ and 
the intersection number  $(C\cdot H)$ is divisible by 
$m\cdot \mtd{n-1}$ for 
some fixed number $m$
for every curve $C\subset X$.
Such examples are  K3 surfaces
(where $m=(H^2)/2$ can take any  value)
or hypersurfaces of very high degree in $\p^4$ \cite{trento-1}.

By Riemann--Roch, $\chi(X,L)-\chi(X,\o_X)$ is divisible by $m$
for any line bundle $L$, thus
$m\ | \ \chi(D, \o_D)$ for every divisor $D\subset X$. 
Thus (\ref{elw.degd.lem}.2)  implies that
$$
\elw_{n-1}(X)\subset \bigl(\elw_{n-2}(X), m\bigr).
$$
 Although
 I do not know if there are further relations between
$\elw_n(X)$ and $\elw_{n-2}(X)$,  the above computations suggest that
once $\elw_n(X)$ and $\elw_{n-2}(X)$
are set, $\elw_{n-1}(X)$ could be any ideal satisfying
$\elw_n(X)\supset\elw_{n-1}(X)\supset\elw_{n-2}(X)$.
\end{exmp}

An extreme case of (\ref{HE.conj}) would be the following.

\begin{prob} Find $n$-dimensional smooth, projective varieties
such that $\elw(X)=1$ yet
 $\elw_i(X)=\bigl(\mtd{n}\bigr)$ for $0\leq i<n$.
\end{prob}

The higher dimensional examples of \cite{trento-1}
 are not convincing for the current purposes, but
the following should be possible to prove.

\begin{prob}  
Given $e_0,e_1, e_2$, 
there is K3 surface $S$ over any field (or even over $\q$) 
 such that $\elw_r(S)=(e_r)$ for $r=0,1,2$ iff
\begin{enumerate}
\item $e_2 \ | \ e_1\ | \ e_0 $
\item  $e_2=\gcd(2, e_1)$
\item $\ e_0\ | \ 12 e_2$ and $\ e_0\ | \ 2 e_1$.
\end{enumerate}
\end{prob}

Another case when (\ref{todd.gives.0.cyc.lem}.1) is known
 connects our index formula
(\ref{rost.prop}) with the version  in \cite{MR3034413}.
The proof given in \cite{MR3034413} relies on \cite{MR1936372};
the claim also follows from properties of the numbers $\mtd{r}$.

\begin{lem} \label{todd.gives.0.cyc.lem.n-1} Let $X$ be a proper, nonsingular
 variety of dimension $n$. Then 
 $$
\mtd{n-1}\cdot \elw_{n-1}(X) \subset \elw_0(X).
$$
\end{lem}

Proof. Let $F$ be a coherent sheaf of dimension $\leq n-1$ on $X$.
Since $X$ is  nonsingular, there are two vector bundles $E_1, E_2$
such that  $[F]=[E_1]-[E_2]$ in $K_0(X)$. Thus, by Riemann--Roch,
$$
\chi(X, F)=\int_X\operatorname{ch}(E_1)\cdot\operatorname{td}(T_X)-
\int_X\operatorname{ch}(E_2)\cdot\operatorname{td}(T_X).
$$
Since $\rank E_1=\rank E_2$, the terms
$\rank E_1\cdot \operatorname{td}_n(T_X)$ and 
$\rank E_2\cdot \operatorname{td}_n(T_X)$ cancel each other.

All other terms have denominators dividing
$(n-r-1)!\cdot \mtd{r}$ for some $r\leq n-1$. 
By  (\ref{mu.td.n}.2) these all divide $\mtd{n-1}$. \qed

%\begin{prob} Connect $\chi(X, \o_X)\in \z/\elw_{\dim X-1}(X)$
%and the Rost invariant.
%Note that the Rost invariant is modulo $\elw_0(X)$.
%It would be very nice to find a coherent interpretation.
%\end{prob}

\section{ELW-indices for special fields}

%\begin{prob} \label{spec.fields.ques}

We consider various 
 fields for which the  sequences
$\elw_0(X),\dots, \elw_{\dim X}(X)$ are special
for all schemes, or at least
for certain interesting classes of varieties.

\subsection*{Finite  fields}{\ }

The following is proved in  \cite{wit}.

\begin{prop}   Let $X$ be a  proper scheme over a finite field $k$. Then
$$
\elw_0(X)= \cdots= \elw_{\dim X}(X)=
\bigl( h^0(\bar Z, \o_{\bar Z}): Z\subset X 
\mbox{ and $Z$ is integral}\bigr).\qed
$$
\end{prop}

The proof is a combination of two lemmas. 

\begin{lem} \label{H0.as.intrscetion.lem}
Let $X$ be a proper $k$-scheme. Then
$$
\elw(X)\subset \bigl( h^0(\bar Z, \o_{\bar Z}): Z\subset X 
\mbox{ and $Z$ is integral}\bigr),
$$
where $Z$ runs through all closed, integral subschemes of $X$ and
$\bar Z\to Z$ denotes the normalization.
\end{lem}

Proof.  Let $W$ be a proper, integral scheme and $F$  a coherent sheaf
 over $W$. Then
$H^i(W, F)$ is a vector space over the field $H^0(W, \o_W) $, 
hence its dimension
is divisible by $h^0(W, \o_W) $. Thus $\chi(W, F)$
is divisible by $h^0(W, \o_W) $. Applying this to $W=\bar Z$
we get that each
$\chi (\bar Z, \o_{\bar Z})$ is contained in
$\bigl( h^0(\bar Z, \o_{\bar Z}): Z\subset X\bigr) $.
By (\ref{deviss.lem}.2) this implies our claim.\qed

\begin{lem}\label{ff.0.cycle}
Let $W$ be a  proper, normal,   integral variety over a finite field $k$.
Then $\elw_0(W)=h^0(W, \o_W)$.  
\end{lem}

Proof. Assume first that $ h^0(W, \o_W)=k$. Then $W$ is 
geometrically integral.
If $\dim W=1$ then by Weil there are points in any large enough field
extension; take two whose degrees are relatively prime.

If $\dim W>1$ then use Bertini to get a geometrically integral
hyperplane section. More precisely, such hyperplane sections exists
over any large enough field
extension; take two whose degrees are relatively prime.

In general, set $K=H^0(W, \o_W)$. After base change to $K$,
the irreducible components  $W_i\subset W_K$ are 
geometrically integral. Thus there is a 0-cycle  $Z_1\subset W_1$
of degree 1. The sum of its conjugates gives a 0-cycle 
of degree $\dim_kK$ on $W$.
 \qed

\subsection*{Henselian fields with algebraically closed residue fields}{\ }

One of the main questions proposed and investigated in
\cite{elw} is the following.

\begin{conj} \label{elw.conj}
 Let $K$ be the quotient field of an excellent,  Henselian DVR with
algebraically closed residue field $k$.
Let $X$ be a proper $K$-scheme. Then 
$$
\elw_0(X)= \elw_1(X)=\cdots= \elw_{\dim X}(X).\eqno{(\ref{elw.conj}.1)}
$$
\end{conj}

The conjecture is almost proved in \cite{elw}. 

\begin{thm} %\cite{elw} 
Let $K$ be the quotient field of an excellent, Henselian DVR with
algebraically closed residue field $k$.
Let $X$ be a proper $K$-scheme. Then
\begin{enumerate}
\item If $\chr k=0$ then 
$\elw_0(X)= \cdots= \elw_{\dim X}(X)$.
\item If $\chr k=p>0$ then 
$\elw_0(X)= \cdots= \elw_{\dim X}(X)$ holds in $\z[p^{-1}]$.
\end{enumerate}
\end{thm}

The key step is the following.

\begin{lem}  %\cite{elw} 
Let $R$ be a  Henselian DVR with quotient field
$K$ and  algebraically closed residue field $k$.
Let $X_R$ be a proper, regular  $R$-scheme with generic fiber $X_K$.
Then $\chi(X_K, \o_{X_K})\in  \elw_0(X_K)$.
\end{lem}

Proof. Write $X_0=\sum_{i\in I} m_i X^i_0$. If $r|m_i$ for every $i$
then  set $Z_r:=\sum (m_i/r) X^i_0$. Note that
$$
\o_X(-iZ_r)^r|_{Z_r}\cong \o_X(-irZ_r)|_{Z_r}\cong\o_{Z_r}.
$$
Thus $\o_X(-iZ_r)|_{Z_r}$ is numerically trivial, hence
 $\chi\bigl(Z_r, \o_X(-iZ_r)|_{Z_r}\bigr)=\chi\bigl(Z_r, \o_{Z_r}\bigr) $. 
Therefore
$$
\chi(X_K, \o_{X_K})=\chi(X_0, \o_{X_0})=
\tsum_{i=1}^r\chi\bigl(Z_r, \o_X(-iZ_r)|_{Z_r}\bigr)=
r\chi\bigl(Z_r, \o_{Z_r}\bigr).
$$
This implies that $\chi(X_K, \o_{X_K})\in (m_i:i\in I)$. 
We conclude by noting that 
through a general point of $X^i_0$ there is a multi-section of degree $m_i$.\qed

\subsection*{Real closed  fields}{\ }

For any scheme $X$ over $\r$, $\elw_0(X)=1$
iff $X(\r)\neq \emptyset$. Otherwise  $\elw_0(X)=2$.
Thus the only question is when the sequence $\elw_i$
drops form 2 to 1.

\begin{exmp} Let $\pi:S\to \p^2$ be a double cover ramified
along a curve of degree $2d$. Let $H$ denote the pull-back of a line in
$\p^2$. 
Then $\pi_*\o_S\cong \o_{\p^2}+\o_{\p^2}(-d)$,
thus
$$
\chi(S, \o_S)=1+\tfrac{(d-1)(d-2)}{2}
\qtq{and} K_S\sim (d-3)H.
$$
If $C\sim rH$ is a curve in $S$ then
$$
\chi(C, \o_C)=r(r+d-3)\tfrac{(H^2)}{2}=r(r+d-3).
$$
Thus if $S(\r)=\emptyset$, $d\equiv 2 \mod 4$ and $\pic(S)=\z[H]$ then
$$
\elw_0(S)=2, \ \elw_1(S)=2,\ \elw_2(S)=1.
$$
Such surfaces can be obtained as small perturbations of
$$
S_0:=\bigl(x^{12}+y^{12}+z^{12}+w^{2}=0\bigr)\subset \p^3(1,1,1,6).
$$
\end{exmp}

Probably there are similar higher dimensional examples.
It is, however, quite difficult to understand all
subvarieties of codimension $\geq 2$ of a given variety.
I do not know how to compute the ELW-indices
for   higher dimensional hypersurfaces.

%The following seems more interesting.

\begin{ques}[OW] \label{RC.over.R.ques}
Let $X$ be a smooth rationally connected variety over $\r$.
Is $\elw_1(X)=1$?
\end{ques}

 It is  known that if $X$ is 
rationally connected and $X(\r)\neq \emptyset$ then it
contains a rational curve \cite[1.7]{MR1715330}. It is conjectured that
$X$ contains a geometrically rational curve even if $X(\r)= \emptyset$;
see \cite[Rem.20]{ar-ko}. (\ref{RC.over.R.ques}) is a weaker
variant of it. 
This is closely related to the conjecture that
the function field of the empty real conic is $C_1$; see
\cite[p.379]{MR0053924}.

%\bibliography{refs-main/refs}
\def\cprime{$'$} \def\cprime{$'$} \def\cprime{$'$} \def\cprime{$'$}
  \def\cprime{$'$} \def\cprime{$'$} \def\dbar{\leavevmode\hbox to
  0pt{\hskip.2ex \accent"16\hss}d} \def\cprime{$'$} \def\cprime{$'$}
  \def\polhk#1{\setbox0=\hbox{#1}{\ooalign{\hidewidth
  \lower1.5ex\hbox{`}\hidewidth\crcr\unhbox0}}} \def\cprime{$'$}
  \def\cprime{$'$} \def\cprime{$'$} \def\cprime{$'$}
  \def\polhk#1{\setbox0=\hbox{#1}{\ooalign{\hidewidth
  \lower1.5ex\hbox{`}\hidewidth\crcr\unhbox0}}} \def\cdprime{$''$}
  \def\cprime{$'$} \def\cprime{$'$} \def\cprime{$'$} \def\cprime{$'$}
\providecommand{\bysame}{\leavevmode\hbox to3em{\hrulefill}\thinspace}
\providecommand{\MR}{\relax\ifhmode\unskip\space\fi MR }
\providecommand{\MRhref}[2]{%
  \href{http://www.ams.org/mathscinet-getitem?mr=#1}{#2}
}
\providecommand{\href}[2]{#2}

\vskip1cm

\noindent Princeton University, Princeton NJ 08544-1000

{\begin{verbatim}kollar@math.princeton.edu\end{verbatim}}

\end{document}